\documentclass[11pt]{article}
\usepackage{amssymb,amsmath,latexsym}
\usepackage{epsfig}

\allowdisplaybreaks

\newtheorem{thm}{Theorem}[section]
\newtheorem{cor}[thm]{Corollary}
\newtheorem{lemma}[thm]{Lemma}

\newtheorem{defn}[thm]{Definition}

\newtheorem{remark}[thm]{Remark}

\numberwithin{equation}{section}

\def\pf{{\medskip\noindent {\bf Proof. }}}
\def\qed{{\hfill $\Box$ \bigskip}}

 \def\sB {{\cal B}} 
 \def\sE {{\cal E}} \def\sF {{\cal F}}
  
  \def\sL {{\cal L}}

 \def\bE {{\mathbb E}}

\def\bP {{\mathbb P}} \def\bQ {{\mathbb Q}} \def\bR {{\mathbb R}}

\def\1{{\bf 1}}

\def\E{{\mathbb E}}
\def\P{{\mathbb P}}

\def\<{\langle}
\def\>{\rangle}
\def\Exp{{\rm Exp}}

\def\K{{\bf K}}
\def\J{{\bf J}}

\def\bea{\begin{align*}}
\def\eea{\end{align*}}
\def\bee{\begin{equation}}
\def\eee{\end{equation}}

\def\eps{\varepsilon}

\def\nin{\noindent }

\makeatletter
\@addtoreset{equation}{section}

\makeatother

\begin{document}
\bibliographystyle{plain}

\title{\Large \bf
Uniform integrability of exponential martingales and
spectral bounds of non-local Feynman-Kac semigroups}

\author{{\bf Zhen-Qing Chen}\thanks{Research partially supported
by NSF Grants DMS-0906743 and DMR-1035196.}}
\date{}
\maketitle

\centerline{Dedicated to Professor Jiaan Yan on the occasion of his 70th birthday}

\bigskip

 \begin{abstract}
In the first part of this paper,  we
give a useful criterion for uniform integrability
 of exponential martingales in the context of Markov processes.
 The condition of this criterion is easy to verify and is, in general,
  much weaker than the commonly used Novikov's condition.
  In the second part of this paper, we present a new approach
  to the study of spectral bounds of Feynman-Kac semigroups
  for a large class of symmetric Markov processes.
  We first establish criteria for the $L^p$-independence
   of spectral bounds for  Feynman-Kac semigroups generated by continuous
   additive functionals, using gaugeability results obtained
   by the author in \cite{C}.
 We then extend these analytic criteria for the
  $L^p$-independence of spectral bounds to
non-local Feynman-Kac semigroups via pure jump Girsanov transforms.
For this, the uniform integrability of the exponential martingales
established in the first part of this paper plays an important role.
We use it to show that Kato classes introduced in \cite{C} can only
 become larger under pure jump Girsanov transforms with symmetric jumping functions.
 \end{abstract}



 \section{Introduction}

Feynman-Kac transform is one of the most important transforms for Markov processes.
Suppose that  $E$ is  a Lusin space
(i.e., a space that is homeomorphic to a Borel subset of a
compact metric space)
and $\sB(E)$ denotes the Borel $\sigma$-algebra on $E$.
Let $m$ be a Borel $\sigma$-finite measure on $E$
with $\hbox{supp} [ m ] =E$ and
$X =(\Omega,\, \sF,\, \sF_t, \, X_t, \, \P_x, \, x\in E)$ be an
$m$-symmetric irreducible
  Borel standard process on $E$ with lifetime $\zeta$
(cf. Sharpe \cite{sharpe} for the terminology). For a continuous additive
functional $A$ of $X$ having finite variations,
one can do Feynman-Kac transform:
$$ T_t f (x)= \E_x \left[ e^{A_t} f(X_t)\right], \qquad t\geq 0.
$$
It is easy to check (see \cite{ABM}) that, under suitable Kato class condition on $A$,
$\{T_t; t\geq 0\}$ forms a strongly continuous symmetric semigroup on $L^p(E; m)$
for every $1\leq p\leq \infty$ and that its $L^2$-infinitesimal generator
is $\sL^\mu:=\sL + \mu$, where $\sL$ is the $L^2$-infinitesimal generator  of the process $X$
and $\mu$ is the (signed) Revuz measure for the continuous additive functional
$A$. To emphasize the correspondence between continuous additive functionals
and Revuz measures, let's denote $A$ by $A^\mu$.
When the process $X$ is discontinuous, it has many discontinuous additive functionals.
Let $F$ be a symmetric function on $E\times E$ that vanishes along the diagonal $d$
of $E\times E$. We always extend it to be zero off $E\times E$.
Then $\sum_{0<s\leq t} F(X_{s-}, X_s)$, whenever it is summable,
 is an additive functional of $X$. Hence one can also perform
 non-local Feynman-Kac transform
 \begin{equation}\label{e:1.1}
  T^{\mu, F}_t f (x):= \E_x \left[ \exp \left( A^\mu_t+\sum_{0<s\leq t} F(X_{s-}, X_s) \right) f(X_t)\right],
   \qquad t\geq 0.
\end{equation}
Non-local Feynman-Kac transforms have been investigated in \cite{C, C2, CS1, CS}.
Let $(N(x, dy), H_t)$ be a L\'evy system of $X$ (see Section \ref{S:2} for its definition).
The infinitesimal generator for $\{T^{\mu, F}_t; t\geq 0\}$ of  \eqref{e:1.1} is
(see Corollary 4.9 and Remark 1 of \cite{CS1})
$$ \sL^{\mu, F}  := \sL + \mu_H {\bf F} + \mu,
$$
where $\mu_H$ is the Revuz measure of the positive continuous additive functional $H$
and
$$\mu_H {\bf F}f (dx):= \left(\int_E \left( e^{F(x, y)}-1\right) f(y) N(x, dy) \right) \mu_H (dx).
$$
Since in this paper we are only concerned with behavior of
the Schr\"odinger semigroup $\{T^{\mu, F}_t; t\geq 0\}$,
by considering the 1-subprocess of $X$ if necessary,  without
loss of generality, we may assume  that $X$ is transient; see Remark \ref{R:4.2} and \eqref{e:5.9} below.
Under some suitable Kato class conditions on the Revuz measure $\mu$ and the function $F$,
$\{T^{\mu, F}_t; t\geq 0\}$ is a strongly continuous symmetric semigroup on
$L^p(E; m)$ for every $1\leq p\leq \infty$. Hence the limit
$$
\lambda_p (X;  \mu+F):= - \lim_{t\to \infty} \frac1t \log \| T^{\mu, F}_t\|_{p, p}
$$
exists, which will be called the {\it $L^p$-spectral bound}
 of the non-local Feynman-Kac semigroup $\{T^{\mu, F}_t; t\geq 0\}$.
We will show in this paper that under suitable  conditions,
 $\lambda_p (X;  \mu+F) = \lambda_2 (X;  \mu+F)$
for all $1\leq p\leq \infty$ if $\lambda_2 (X;  \mu+F)\leq 0$.
If in addition $X$ is conservative, then $\lambda_2 (X;  \mu+F)\leq 0$
becomes a necessary and sufficient condition for the independence of
$\lambda_p (X, \mu+F)$ in $p\in [1, \infty]$.
The $L^2$-spectral bound $\lambda_2 (X;  \mu+F)$ has a
 variational formula in terms of the Dirichlet form of $X$, $\mu$ and $F$, see
\eqref{e:5.7} below.
 The spectral bound results obtained in this paper not only   extend
   earlier results in \cite{T, T2, TT, TT3, Taw} to
   a larger class of symmetric Markov processes
   but also give several new criteria (for example, Theorem
    \ref{T:4.7}(i),
   Theorem \ref{T:4.8} and Theorems \ref{T:5.3}-\ref{T:5.5}).
   See Remarks \ref{R:4.10}  and \ref{R:5.6} below for details.

When $F=0$, the $L^p$-independence of spectral bounds
for continuous Feynman-Kac transforms was investigated by
Takeda in \cite{T, T2} for conservative Feller processes and for symmetric
Markov processes with strong Feller property and a tightness assumption,
 respectively,
using a large deviation approach. The results in \cite{T} were extended
to purely discontinuous Feynman-Kac transforms (i.e. with $\mu=0$) first
in \cite{TT} for rotationally symmetric $\alpha$-stable processes
and then in \cite{Taw} for   conservative
doubly Feller processes,
both papers  again using a large deviation approach.
A stochastic process is said to be doubly Feller if it is a
Feller process having strong Feller property.
See also \cite{DKK} for further extensions of above results
for doubly Feller processes, which
are established along a similar line using large deviation approach.
The approach of this paper is completely different.
We use the   gaugeability results obtained in \cite{C} to establish
the $L^p$-independence of spectral bounds for local
Feynman-Kac semigroups for a large class of symmetric Markov processes.
These results extend the main results in \cite{T, T2}.
We then show that using a pure jump
Girsanov transform, we can reduce a non-local Feynman-Kac
transform for $X$ into a
continuous Feynman-Kac transform for the Girsanov transformed process $Y$
and then apply the $L^p$-independence result for local Feynman-Kac
semigroups. For this,
uniform integrability of the exponential martingale used in the
Girsanov transform plays a crucial role. Thus in the first part of this paper,
we present a  useful criterion for the uniform integrability
 of exponential martingales in the context of Markov processes, which is of
 independent interest.
 The condition of this criterion is easy to verify and is, in general,
  much weaker than the commonly used Novikov's condition.
  The special cases of this criterion have been used earlier
  in \cite{CWZ} and \cite{C2}.
  Using a super gauge theorem established in \cite{C}, we show that
  the Kato classes of $X$  introduced in \cite{C} are contained in the corresponding
  Kato classes of the Girsanov transformed process $Y$.

  The rest of the paper is organized as follows.
  In Section \ref{S:2}, we give precise setup of this paper, including the
  definitions of Kato classes and L\'evy systems.
  The criterion of the  uniform integrability
 of exponential martingales in the context of Markov processes is presented
 and proved in Section \ref{S:3}. Spectral bounds for local
 Feynman-Kac semigroups and its $L^p$-independence are studied in Section \ref{S:4}, using
 gaugeability results for Feynman-Kac transforms obtained by the author in \cite{C}.
In Section \ref{S:5}, we first show that the Kato classes of $X$
 are contained in the corresponding Kato classes of the Girsanov transformed process $Y$,
 and then use it derive the criteria for the $L^p$-independence of spectral bounds
 for non-local Feynman-Kac semigroups.
 To keep the exposition of this paper as transparent as possible,
 we have not attempted to present the most general conditions
 on $\mu$ and $F$.

\section{Kato classes and non-local Feynman-Kac transform}\label{S:2}

Let $E$ be a Lusin space
and $\sB(E)$ be the Borel $\sigma$-algebra on $E$.
Let $m$ be a Borel $\sigma$-finite measure on $E$
with $\hbox{supp} [ m ] =E$ and
$X =(\Omega,\, \sF,\, \sF_t, \, X_t, \, \P_x, \, x\in E)$ be an
$m$-symmetric irreducible
transient Borel standard process on $E$ with lifetime $\zeta$.
As mentioned in the Introduction, the transience assumption on $X$
here is just for convenience and is unimportant---we can always consider the 1-subprocess of $X$
instead of $X$ if necessary; see Remark \ref{R:4.2} and \eqref{e:5.9}.
Let $(\sE, \sF)$ denote the Dirichlet form of $X$;
that is, if we use $\sL$ to denote the infinitesimal generator of $X$,
then $\sF$ is the domain of the operator $\sqrt{-\sL}$ and
for $u, v \in \sF$,
$$ \sE (u, v) = (\sqrt{-\sL u}, \, \sqrt{-\sL v} )_{L^2(E, m)}.
$$
We refer readers to \cite{CF} or \cite{FOT} for terminology and various
properties of Dirichlet forms such as continuous additive functional,
martingale additive functional.

The transition operators $\{P_t, t\ge 0\}$ of $X$ are defined by
$$
P_tf(x):=\E_x [f(X_t)] =\E_x[ f(X_t); \, t<\zeta ].
$$
(Here and in the sequel, unless mentioned otherwise,
we use the
convention that a function defined on $E$ takes the value 0 at the
cemetery point $\partial$.)
Throughout this paper,
 we assume that there is a Borel symmetric function $G(x, y)$
on $ E\times E$ such that
$$
\E_x \left[ \int_0^\infty f(X_s) ds \right]
=\int_E G(x, y) f(y) m(dy)
$$
for all measurable function $f\geq 0$.
$G(x, y)$ is called the Green function of
$X$. The Green function $G$ will always be chosen so that for each fixed
$y\in E$,  $x\mapsto G(x, y)$ is an
excessive function of $X$.
Note that we do not assume $X$ is a Feller process, nor do we assume
$X$ has strong Feller property.

For every $\alpha>0$, one deduces from
the existence of the Green function $G(x, y)$
that there exists
a kernel $G_\alpha(x, y)$ so that
$$
\E_x \left[ \int_0^\infty e^{-\alpha s} f(X_s) ds \right]
=\int_E G_\alpha (x, y) f(y) m(dy)
$$
for all measurable $f\geq 0$. Clearly, $G_\alpha (x, y)\leq G(x, y)$.
Note that by \cite[Theorem 4.2.4]{FOT}, for every $x\in E$ and $t>0$,
$X_t$ under $\P_x$ has a density function $p(t, x, y)$
with respect to the measure $m$.

A set $B$ is said to be $m$-polar if
$\P_m (\sigma_B<\infty )=0$, where
$\sigma_B:=\inf \{ t>0: \, X_t \in B \}$.
We call a positive measure $\mu$ on $E$ a smooth measure
of $X$
if there is a positive continuous additive functional
(PCAF in abbreviation)
$A$ of $X$ such that
\begin{equation}\label{eqn:Revuz1}
\int_E f(x) \mu (dx) = \uparrow \lim_{t\downarrow 0}
\E_m \left[ \frac1t \int_0^t f(X_s) dA_s \right].
\end{equation}
for any Borel $f\geq 0$.
Here $\uparrow \lim_{t\downarrow 0}$ means the quantity
is increasing
as $t\downarrow 0$.
The measure $\mu$ is called the
Revuz measure of $A$.
We refer to \cite{CF, FOT}
 for the characterization of smooth
measures in terms of nests and capacity.

For any given positive smooth measure $\mu$, define
$G \mu (x) =\int_E G(x, y) \mu (dy)$.
It is known (see Stollmann and Voigt \cite{SV})
that for any positive smooth measure
$\mu$ of $X$,
\begin{equation}\label{eqn:en}
\int_E u(x)^2 \mu(dx) \leq \|G\mu \|_\infty \,
\sE (u, u) \qquad \hbox{for } u \in \sF.
\end{equation}
Recall that as $X$ is assumed to have a Green function,
any $m$-polar set is polar. Hence by
(\ref{eqn:en}) a PCAF $A$ in the sense of
\cite{FOT} with an exceptional set that has a bounded potential
(that is, $x \mapsto \E_x \left[ A_\zeta \right]=G \mu$ is bounded
almost everywhere on $E$,
where $\mu$ is the Revuz measure of $A$)
can be uniquely
refined into a PCAF in the strict sense (as defined on p.195 of
\cite{FOT}).
This can be proved by using the same argument as that in the
 proof of Theorem 5.1.6 of   \cite{FOT}.

\medskip

The following definitions are taken from Chen \cite{C}.

\begin{defn}\label{D:2.1}
Suppose that $\mu$ is a signed smooth measure.
Let $A^\mu$ and $A^{|\mu|}$ be the continuous additive
functional and positive continuous additive functional
  of $X$ with Revuz measures
 $\mu$ and $|\mu |$, respectively.

\begin{description}
\item{\rm (i)} We say  $\mu$ is
in the Kato class of $X$, $\K(X)$ in abbreviation,
if
$$ \lim_{t\to 0} \sup_{x\in E}
\E_x \left[ A^{|\mu |}_t \right] =0.
$$

\item{\rm (ii)} $\mu$  is said to be in the class $\K_\infty (X)$
if for any $\eps>0$, there is a Borel set
$K=K(\eps)$ of finite $|\mu|$-measure
and a constant $\delta=\delta(\eps)>0$
such that for all measurable set $B\subset K$
with $|\mu|(B)<\delta$,
\begin{equation}\label{eqn:1.1}
\| G (\1_{K^c\cup B} | \mu |) \|_\infty <\eps.
\end{equation}

\item{\rm (iii)} $\mu$ is said to be in the class $\K_1 (X)$
if there is a Borel set $K$ of finite $|\mu|$-measure
and a constant $\delta>0$
such that
\begin{equation}\label{eqn:1.3}
\beta_1 (\mu):= \sup_{B\subset K: \, |\mu|(B)<\delta}
\| G (1_{K^c \cup B} | \mu |) \|_\infty <1.
\end{equation}

\item{\rm (iv)} A function $q$ is said to be in class
$\K (X)$,  $\K_\infty (X)$ or $\K_1 (X)$
 if $\mu(dx) :=q(x) m(dx)$ is in the corresponding spaces.
\end{description}
\end{defn}

According to \cite[Proposition 2.3(i)]{C}, $\K_\infty (X)\subset \K (X)\cap \K_1(X)$.
Suppose that $\mu$ is a positive measure in $\K_1 (X)$.
By Propositions 2.2 in \cite{C}, $G\mu (x)=\E_x [ A^\mu_\infty]$
 is bounded
and so (\ref{eqn:en}) is satisfied. Therefore the PCAF
corresponding to $\mu$ can and is always taken to be
in the strict sense.

\medskip

Let $(N, H)$ be a L\'evy system for $X$
(cf. Benveniste and Jacod \cite{BJ}
and Theorem 47.10 of Sharpe \cite{sharpe});
that is,
$N(x, dy)$  is a kernel from $(E , {\cal B}(E))$ to
$(E, {\cal B}(E))$ satisfying
$N(x, \{ x \})=0$,
and $H_t$ is a PCAF of $X$ with bounded
1-potential such that  for any nonnegative Borel function $f$ on $E \times
E$ vanishing on the diagonal and any $x\in E$,
\begin{equation}\label{eqn:LS1}
\E_x\left[ \sum_{s\le t}f(X_{s-}, X_s) \1_{\{s<\zeta\}}\right]
=\E_x\left[\int^t_0\int_E f(X_s, y) N(X_s, dy)dH_s\right].
\end{equation}
The Revuz measure for $H$ will be denoted as $\mu_H$.

\medskip

\begin{defn}\label{D:2.2}
Suppose $F$ is a bounded function
on $E\times E$ vanishing on the diagonal $d$.
It is always extended to be zero off $E\times E$.
Define $\mu_F (dx):= \left(
\int_E F(x, y) N(x, dy)\right) \mu_H (dx)$.
We say $F$ belongs to the class $\J(X)$ (respectively, $\J_{\infty}(X)$) if
the measure
$$\mu_{|F|} (dx):= \left( \int_{E} |F(x, y)| N(x, dy)\right)
                     \mu_H (dx)
$$
belongs to $\K(X)$ (respectively, $\K_{\infty}(X)$).
\end{defn}

\medskip

See \cite[Section 2]{C2} for concrete examples of
$\mu \in \K_\infty (X)$ and $F\in \J_\infty (X)$.
The following result is established in \cite[Theorems 2.13 and 2.17]{C}.

\begin{thm}\label{T:2.3}
Assume that a signed measure
$\mu \in \K_\infty (X)$ and $F\in \J_{\infty}(X)$.
Let $A^\mu$ be the continuous additive functional of $X$
with signed Revuz measure $\mu$, and define
the non-local Feynman-Kac functional
$$
  e_{A^\mu+F}(t):= \exp \Big( A^\mu_t + \sum_{0<s\leq t}
  F(X_{s-}, X_s)\Big),   \qquad t\geq 0 .
$$

\begin{description}
\item{\rm (i)} {\rm (Gauge Theorem)}
The gauge function $g(x):=\E_x \left[ e_{A^\mu+ F} (\zeta) \right]$
is either bounded on $E$ or identically $\infty $ on $E$.
When $g$ is bounded,  we say $(X, A^\mu+F)$,  or $(X, \mu+F)$,
  is gaugeable.

\item{\rm (ii)} {\rm (Super Gauge Theorem)}
Suppose that $(X, A^\mu+F)$ is gaugeable.
Then there is an $\eps_0>0$ such that
$(X, A^{\mu+\eps_0 |\mu|}+ F+\eps_0|F|)$
is gaugeable. In particular,
$(X, A^{(1+\eps)\mu}+(1+\eps)F)$ is gaugeable
for all $\eps\in [0, \eps_0]$.
\end{description}
\end{thm}

\section{Uniform integrability of exponential martingales}\label{S:3}

Our uniform integrability criterion for Dol\'eans-Dade exponential
martingales is based on the following simple observation.

\begin{lemma}\label{L:3.1}
 Suppose that $Z=\{Z_t, \sF_t\}_{t\geq 0}$ is a non-negative supermartingale and define $Z_\infty=\lim_{t\to \infty} Z_t$.
If there is a constant $c>0$ so that $Z_t\leq c\, \E [Z_\infty |\sF_t]$,
then $Z$ is uniformly integrable.
\end{lemma}

\pf Since $Z$ is a non-negative supermartingale, $Z_\infty=\lim_{t\to
\infty} Z_t$ exists a.s. and so, by Fatou's lemma,
$\E [ Z_\infty] \leq \E [ Z_0]<\infty$. The conclusion of the lemma follows
from the fact that $\{ \E [Z_\infty | \sF_t], t\geq 0\}$ is uniformly integrable. \qed

Suppose that $M=\{M_t, \sF_t\}_{t\geq 0}$ is a local martingale with $M_0=0$.
It is well-known (see, e.g., \cite[Theorem 9.39]{HWY}) that
$$ Z_t=1+ \int_0^t Z_{s-} dM_s
$$
has a unique solution, which is given by
\begin{equation}\label{e:3.1}
Z_t= \exp \left(M_t-\frac12 \< M^c\>_t\right) \prod_{0<s\leq t} (1+\Delta M_s) e^{-\Delta M_s}
\qquad \hbox{for } t\geq 0.
\end{equation}
Here $M^c$ is the continuous local martingale part of $M$ and $\<M^c\>$ is
the quadratic variation  process of $M^c$.
The quadratic variation process $[M]$ of $M$ is defined as
$$ [M]_t=\<M^c\>_t+ \sum_{0<s\leq t} (\Delta M_s)^2, \qquad t\geq 0.
$$
The local martingale $Z$ is called the Dol\'eans-Dade exponential martingale of $M$
and will be denoted as $\Exp (M)$.

In the remainder of this section, $X$ is a general strong Markov process
 (not necessarily symmetric).

\begin{thm}\label{T:3.2}
 Suppose that $M$ is a martingale additive functional
 of a strong Markov process $X$ with $M_0=0$ and
$\sup_{x\in E} \E_x [ M]_\infty  <\infty$ and that there is
a constant $\delta \in (0, 1)$ so that $\Delta M_s \geq \delta -1$
for every $s\geq 0$ a.s.
Then $\Exp (M)$ is uniformly integrable
under $\P_x$ for every $x\in E$.
\end{thm}

\pf Under the assumption that $\sup_{x\in E}  \E_x [ M]_\infty <\infty$,
$M$ is a square-integrable martingale under $\P_x$ for every $x\in E$.
Thus $M_\infty = \lim_{t\to \infty} M_t$ exists $\P_x$-a.s. and $\E_x [M_\infty]=0$
for every $x\in E$.
Observe that $Z:=\Exp (M)$ is a non-negative local martingale and hence a non-negative
supermartingale under $\P_x$ for every $x\in E$.
It is also a multiplicative functional of $X$.
 By \eqref{e:3.1}, the Markov
property of $X$ and Jensen's inequality, for every $x\in E$ and $t\geq 0$,
\begin{eqnarray*}
 \E_x \left[ Z_\infty/ Z_t \, | \, \sF_t\right]
  &=& \E_x \left [Z_\infty \cdot \theta_t | \sF_t\right] = \E_{X_t} [Z_\infty] \\
&\geq& \exp \Big( \E_{X_t}\Big[ M_\infty -\frac12 \<M^c\>_\infty
 + \sum_{0<s<\infty} (\log (1+\Delta M_s) -\Delta M_s )\Big] \Big)\\
 &\geq& \exp \Big( \E_{X_t}\Big[ -\frac12 \<M^c\>_\infty
 -\frac{1}{2\delta^2} (\Delta M_s )^2\Big] \Big)\\
 &\geq &  \exp \Big(   -\frac1{2\delta^2} \sup_{x\in \infty} \E_x [ M ]_\infty \Big),
\end{eqnarray*}
where the second to last inequality is due to the fact that
 $$ \log (1+x)-x \geq x^2/(2\delta^2) \qquad \hbox{for } x\geq -1+\delta.
 $$
The conclusion of the theorem now follows from Lemma \ref{L:3.1}. \qed

\begin{cor}\label{C:3.3}
 Suppose that  $X$ is a
 transient
 Borel standard process on $E$ with L\'vey system $(N, H)$.
  Let $b$ be a function
on $E  \times E $ vanishing on the diagonal $d$ such that
$b(x, y)\geq  \delta -1$ for some constant $\delta>0$ and that
$G\mu_{b^2}$ is bounded for
$$ \mu_{b^2}(dx):= \left( \int_E  b (x, y)^2 N(x, dy)\right) \mu_H(dx).
$$
 Then there is a unique purely discontinuous
square integrable martingale additive functional of $X$ with
$M_0=0$ and $\Delta M_t = b(X_{t-}, X_t)\, \1_{\{t<\zeta\}}$ for $t>0$. Moreover,
$\Exp (M)$ is uniformly integrable
under $\P_x$ for every $x\in E$.
\end{cor}

\pf Using the L\'evy system, the assumption that $\mu_{b^2}$ has a bounded
potential is equivalent to the assumption that
$$ \sup_{x\in E} \E_x \left[ \sum_{s>0} b(X_{s-}, X_s)^2
\, \1_{\{s<\zeta\}}\right] <\infty.
$$
Thus there is a unique purely discontinuous
square integrable martingale additive functional of $X$ with
$M_0=0$ and $\Delta M_t = b(X_{t-}, X_t) \, \1_{\{t<\zeta\}}$
for $t>0$ (see, e.g., \cite{HWY}).
Since
$$
[M]_\infty =\sum_{s>0} b(X_{s-}, X_s)^2 \, \1_{\{s<\zeta\}},
$$
it follows immediately from Theorem \ref{T:3.2} that $\Exp (M)$
is uniformly integrable. \qed

\begin{remark}\label{R:3.4} \rm
(i) Suppose that $E=D$ is a connected open subset of $\bR^n$, $X$
is Brownian motion killed upon leaving domain $D$, and $M_t=\int_0^t b(X_s) dX_s$.
Note that $M_t=\int_0^{t \wedge \tau_D}b(X_s) dX_s$ is a martingale additive
functional of $X$, where $\tau_D$ is the first exit time (or lifetime) of
the killed Brownian motion $X$ from $D$.
We have by Theorem \ref{T:3.2} that
  the exponential martingale $\Exp (M)$ is uniformly integrable
  if
  $$\sup_{x\in D} \E_x \left[ \int_0^\infty |b(X_s)|^2 ds\right]
  = \sup_{x\in D} \E_x \left[ \int_0^{\tau_D} |b(X_s)|^2 ds\right]<\infty.
  $$
The result in this particular case
was first derived in passing on page 746 of \cite{CWZ}.
The above condition is in general much weaker than
the Novikov's condition for the
uniform integrability of $\Exp (M)$.
When $M$ is a  continuous local martingale $M$, for $\Exp (M)$ to
be uniformly integrable, Novikov's condition requires
$\E_x \exp (\< M\>_\infty/2)<\infty$ (see, e.g.,
\cite[Proposition VIII.1.15]{RY}).
In this concrete example, the latter
condition  amounts to the assumption that
$$
\E_x \exp \Big(\frac12 \int_0^\infty |b(X_s)|^2 ds\Big)
=\E_x \exp \Big(\frac12 \int_0^{\tau_D} |b(X_s)|^2 ds\Big)  <\infty.
$$

(ii) Suppose that  $X$ is an $m$-symmetric irreducible transient Borel standard process on $E$, and
$F\in \J(X)$.
Then
\bee \label{e:3.2}
 M_t:= \sum_{0<s\leq t}
 F(X_{s-}, X_s) - \int_0^t \left( \int_{E_\partial} F(X_s, y)
 N(X_s, dy) \right)dH_s, \qquad t\geq 0,
\eee
is the purely discontinuous martingale additive functional of $X$ with $M_0=0$ and
$\Delta M_t=b(X_{t-}, X_t)$ for $t>0$. Define
\bee \label{e:3.3}
A^F_t:= \int_0^t \left( \int_E \Big( e^{F(x, y)}-1 \Big) N(X_s, dy)\right) dH_s,
\eee
which is a continuous additive functional of $X$.
Then by  \eqref{e:3.1}, we have
\bee \label{e:3.4}
 \Exp (M)_t =\exp \Big( \sum_{s\leq t}  F(X_{s-}, X_s)- A^F_t \Big),
 \qquad t\geq 0.
\eee
Suppose now that $F\in \J_\infty (X)$. Then
  the condition of Corollary \ref{C:3.3}
is satisfied for $b(x, y):=e^{F(x, y)}-1$ and so
 $\Exp (M)$ is uniformly integrable under $\P_x$
for every $x\in E$.
This fact was first established in \cite[page 241]{C2}.

(iii) If the condition $\sup_{x\in E} \E_x [ M]_\infty <\infty$
in Theorem \ref{T:2.3} is replaced by a weaker condition
 $\sup_{x\in E} \E_x [ M]_T <\infty$ for some fixed constant $T>0$,
the same proof of Theorem \ref{T:2.3} yields that
$\{ \Exp (M)_t, t\in [0, T]\}$ is a $\P_x$-martingale for every
$x\in E$.    \qed
\end{remark}

The above uniform integrability results for exponential martingales
 will be used in Section \ref{S:5}.

\section{Spectral bounds for local Feynman-Kac semigroup}\label{S:4}

For a signed measure $\mu$, we use $\mu^+$ and $\mu^-$ to denote
the positive part and negative part of $\mu$  appearing in
the Hahn-Jordan decomposition of $\mu$.
Observe that if $\mu_1$ and $\mu_2$ are two non-negative measures
so that $\mu_1-\mu_2=\mu$, then $\mu_1 \geq \mu^+$ and
$\mu_2\geq \mu^-$.

In the rest of this paper, we work under the setting of Section 2.
Let $\mu$  be a signed smooth measure so that $\mu^+ \in \K_1 (X)$
and $G\mu^-$ is bounded.
Define the Feynman-Kac semigroup $\{T^\mu_t, t\geq 0\}$ by
$$ T^\mu_t f(x)=\E_x \left[ e^{A^\mu_t} f(X_t) \right] .
$$
As is explained in the paragraph proceeding \cite[Theorem 2.12]{C},
 $\{T^\mu_t, t\geq 0\}$   is a
strongly continuous symmetric semigroup in $L^p(E, m)$ for every $1\leq p\leq \infty$ and its associated symmetric
quadratic form is $(\sE^\mu, \sF)$, where
$$ \sE^\mu (u, v)= \sE(u, v) -\int_E u(x) v(x) \mu (dx)
\qquad \hbox{for } u, v\in \sF.
$$
If we use ${\cal L}$ to denote the infinitesimal generator for the
semigroup of $X$, then the infinitesimal generator for the semigroup
$P^\mu_t$ is ${\cal L}+\mu$.
For $1\leq p\leq 1$, we use $\| T^\mu_t\|_{p, p}$ to denote
the operator norm of $T^\mu_t: L^p (E; m)\to L^p (E; m)$.

\bigskip

For $1\leq p \leq \infty$, define the $L^p$-spectral
bound of semigroup $\{T^\mu_t, t\geq 0\}$ by
$$ \lambda_p(X, \mu):= - \lim_{t\to \infty} \frac1t
\log \|T^\mu_t \|_{p,p} = - \inf_{t>0} \frac1t \log
\|T^\mu_t \|_{p,p}.
$$
Clearly
$$ \| T^\mu_t \|_{\infty, \infty} =\| T^\mu_t 1\|_\infty =\sup_{x\in E}
\E_x \left[ e_A(t); t<\zeta  \right].
$$
It is well-known  that the $L^2$-spectral bound $\lambda_2(X, \mu)$
of $\{T^\mu_t, t\geq 0\}$ can be represented in terms of its quadratic form
$(\sE^\mu, \sF)$:
\begin{eqnarray}\label{e:4.1}
\lambda_2 (X, \mu )&=& \inf\left\{ \sE^\mu (u, u): \,
      u\in \sF \ \mbox{ with } \int_E u(x)^2 m(dx)=1 \right\} \nonumber\\
&=&\inf \left\{ \sE(u, u)-\int_E u(x)^2 \mu (dx):
\, u\in \sF  \ \mbox{ with } \int_E u(x)^2 m(dx)=1 \right\}.
\end{eqnarray}
By duality, we have  $\| T^\mu_t\|_{1,1}= \| T^\mu_t \|_{\infty, \infty}$.
Consequently, it follows from the Cauchy-Schwartz inequality that
$$
 \|T^\mu_t f\|^2_2 \leq  \| T^\mu_t 1\|_\infty \, \| T^\mu_t (f^2)\|_1
\leq \| T^\mu_t\|_{\infty, \infty}^2 \| f\|^2_2 \qquad \hbox{for } f\in L^2(E, m).
$$
Thus we have $\| T^\mu_t \|_{2, 2}  \leq \|T^\mu_t\|_{\infty, \infty}$.
We now deduce by interpolation that
$$
\| T^\mu_t \|_{2,2}\leq \| T^\mu_t\|_{p,p} \leq \| T^\mu_t \|_{\infty, \infty}
\quad \mbox{for } 1< p< \infty.
$$
Hence
\begin{equation}\label{e:4.2}
\lambda_\infty (X, \mu) \leq \lambda_p (X, \mu)
\leq \lambda_2 (X, \mu) \quad \mbox{for } 1<p<\infty.
\end{equation}

\medskip

The following theorem is   proved as Theorem 2.12 in  \cite{C}, however
  condition \eqref{e:lp3} is missing from
its statement.
 For reader's convenience, we reproduce the proof here.

\begin{thm}\label{T:eigen}
Assume that   $m(E)<\infty$,  $\|G1\|_\infty <\infty$.
Let $\mu$ be  a signed smooth measure
such that  $\mu^+ \in \K_1 (X)$ and $G\mu^-$ is bounded.
Then $(X, \mu)$ is gaugeable   if and only if $\lambda_2(X, \mu)>0$.
Assume in addition that $\mu^+\in \K(X)$ and   that
\begin{equation}\label{e:lp3}
 \hbox{there is some } t_0> 0 \hbox{ so that }
 P_{t_0} \hbox{ is a bounded operator from }
 L^2(E; m) \hbox{ into } L^\infty (E; m).
\end{equation}
Then $\lambda_p(X, \mu)$ is independent of $p\in [1, \infty]$
 if $\lambda_2(X, \mu)>0$.
\end{thm}

\pf
By \cite[Theorem  2.11]{C}, if $(X, \mu)$ is gaugeable,
then $\lambda_\infty (X, \mu)>0$ and therefore
$\lambda_2 (X, \mu)>0$.
Conversely suppose  $\lambda_2 (X, \mu)>0$. Then
for any $\eps \in (0, \lambda_2(X, \mu))$,
there is $\delta (\eps) >0$ such that
\begin{equation}\label{eqn:eig2}
 \| T^\mu_t \|_{2, 2} \leq e^{-t(\lambda_2(X, \mu) -\eps )} \quad
\mbox{for } t\geq \delta (\eps).
\end{equation}
Since $1\in L^2(E, m)$,  $\int_0^\infty T^\mu_t 1  \, dt$
is $L^2(E; m)$-integrable. Hence by \cite[Theorem 2.11]{C},
$(X, \mu)$ is gaugeable.

Assume now  that \eqref{e:lp3} holds and that $\mu^+\in \K(X)$.
By duality, $P_{t_0}$ is a bounded operator from $L^1(E; m)$ to $L^2(E; m)$.
It follows that $P_{2t_0}=P_{t_0}\circ P_{t_0}$ is a bounded operator
from $L^1(E; m)$ into $L^\infty (E; m)$, whose operator norm will be
denoted as $\| P_{2t_0}\|_{1, \infty}$. On the other hand, since $\mu^+\in
\K (X)$, there is some $\delta >0$ so that
$ \sup_{x\in E}\E_x [ A^{\mu^+}_\delta]<1/2$. By Khasminskii's inequality,
$$
c_1:= \sup_{x\in E} \E_x \left[ \exp (2A^{\mu^+}_\delta) \right]
\leq \frac{1}{1- \sup_{x\in E}\E_x \big[ 2A^{\mu^+}_\delta \big]}
<\infty.
$$
Thus by the Markov property of $X$, we have
$\sup_{x\in E} \E_x \left[ \exp (A^{\mu^+}_t) \right]<\infty$
 for every $t>0$. Thus for every $f\in L^2(E; m)$ and $x\in E$,
 by Cauchy-Schwartz inequality,
 \begin{equation}\label{e:4.5a}
 |T^\mu_{2t_0} f(x)|^2
 =\left( \E_x \left[ \exp (A^\mu_{2t_0}) f(X_{2t_0}) \right] \right)^2
 \leq \E_x \left[ f(X_{2t_0})^2\right] \,
 \E_x \left[ \exp (2A^{\mu^+}_{2t_0}) \right]
 \leq c_1 \, \|P_{2t_0}\|_{1, \infty} \, \|f \|_2^2.
 \end{equation}
Suppose $\lambda_2(X, \mu)>0$. For any
$\eps \in (0, \lambda_2(X, \mu))$, there is $\delta (\eps)>0$  so that
\eqref{eqn:eig2} holds. Then for  $t>\delta(\eps)+2t_0$, by
  \eqref{e:4.5a} and then \eqref{eqn:eig2},
\begin{eqnarray*}
 \| T^\mu_t \|_{\infty , \infty}&=&\|T^\mu_t 1\|_\infty =
 \| T^\mu_{2t_0} (T^\mu_{t-2t_0} 1)\|_\infty
  \leq c_2 \, \| T^\mu_{t-2t_0} 1)  \|_2 \leq c_2\,  \sqrt{ \, m(E)} \,
e^{-(t-1)(\lambda_2(X, \mu)-\eps)}.
\end{eqnarray*}
This implies that $\lambda_\infty (X, \mu)\geq \lambda_2(X, \mu)-\eps$
and so $\lambda_\infty (X, \mu) \geq \lambda_2(X, \mu)$.
Hence by (\ref{e:4.2}),
$\lambda_\infty(X, \mu)=\lambda_2(X, \mu)=\lambda_p(X, \mu)$
for all $p\in [1, \infty]$.
\qed

\bigskip

For $\alpha>0$, let $X^{(\alpha)}$ denote the $\alpha$-subprocess of $X$;
that is, $X^{(\alpha)}$ is the  subprocess of $X$
killed at exponential rate $\alpha$.
Let $G^{(\alpha)}$ be the $0$-resolvent (or Green operator) of $X^{(\alpha)}$.
Then $G^{(\alpha)}=G_\alpha$,  the $\alpha$-resolvent of $X$.
Thus for $\beta >\alpha >0$,
 $ \K_1 (X)\subset \K_1 (X^{(\alpha )}) \subset
\K_1 (X^{(\beta )})$ and
$\K_\infty (X)\subset \K_\infty (X^{(\alpha )}) \subset
\K_\infty (X^{(\beta )}) $.
In fact, it follows from the resolvent equation
$G_\alpha =G_\beta + (\beta -\alpha) G_\alpha G_\beta$
that $\K_\infty (X^{(\alpha)}) = \K_\infty (X^{(\beta)})$
for every $\beta >\alpha$.
Consequently,  $\J_\infty (X^{(\alpha)}) = \J_\infty (X^{(\beta)})$
for every $\beta >\alpha$.

\medskip

\begin{remark}\label{R:4.2}\rm
For any signed measure $\mu$ on $E$ with
$\mu \in \K_1 (X^{(\alpha)})$ and $G_\alpha \mu^-$ bounded
for some $\alpha >0$,
$\{T^\mu_t, t\geq 0\}$ is still well defined as a strongly
continuous symmetric
semigroup in $L^p(E; m)$ for every $p\in [1, \infty]$,
and relation \eqref{e:4.2} continues to hold.
To see this,
let $\{T^{\mu, \alpha}_t, t\geq 0\}$ denote the Feynman-Kac semigroup
of $X^{(\alpha)}$ associated with smooth measure $\mu$.
Recall the following simple facts.
Let $R$ be an exponential random variable with mean $1/\alpha$ that is independent of $X$.
The $\alpha$-subprocess $X^{(\alpha)}$ of $X$ can be realized as follows:
$X^{(\alpha)}_t (\omega)= X_t (\omega)$ if $t<R(\omega)$ and
$X^{(\alpha)}_t =\partial$ if $t\geq R(\omega)$.
Let $A^\mu$ be the continuous additive functional of $X$ with (signed)
Revuz measure $\mu$. Then $t\mapsto A_{t\wedge R}$ is a continuous
additive functional of $X^{(\alpha)}$ with Revuz measure $\mu$.
It follows immediate that $\K(X)=\K(X^{(\alpha)})$ and
$$ T^{\mu, \alpha}_t =e^{-\alpha t} T^\mu_t
\qquad \hbox{for every } t\geq 0.
$$
Since $\{T^{\mu, \alpha}_t, t\geq 0\}$  is a strongly continuous symmetric semigroup in $L^p(E; m)$ for every $1\leq p\leq \infty$, so is $\{T^\mu_t, t\geq 0\}$.
Moreover,
\begin{equation}\label{e:radius}
\lambda_p (X^{(\alpha)}, \mu) = \alpha
 + \lambda_p (X, \mu) \qquad \hbox{for every }
p\in [1, \infty] \hbox{ and } \beta >0.
\end{equation}
Since relation \eqref{e:4.2} holds for $\lambda_p (X^{(\alpha)},
\mu)$, the same holds for $\lambda_p (X, \mu)$.
\end{remark}

\begin{cor}\label{C:radius}
Assume that \eqref{e:lp3} holds and that $m(E)<\infty$.
Let $\mu$ be  a signed smooth measure
with $\mu^+\in \cup_{\alpha>0}\K_1 (X^{(\alpha)}) \cap \K (X)$
and $G_{\alpha_0} \mu^-$ bounded for some $\alpha_0>0$.
Then
$\lambda_p(X, \mu)$ is independent of $p\in [1, \infty]$.
\end{cor}

\pf There is an $\alpha >0$ sufficiently large so that $\mu^+ \in \K_1 (X^{(\alpha)})\cap \K (X)$ with $G_\alpha \mu^-$ bounded.
Note that since $\mu \in \K (X )$,
 $\lambda_2 (X, \mu)$ is finite. By increasing the value of $\alpha$ if necessary,
 we may and do assume $\lambda_2 (X^{(\alpha)}, \mu )>0$.
 On the other hand, $G^{(\alpha)}1 \leq 1/\alpha$. Therefore by Theorem \ref{T:eigen},
 we have $\lambda_2 (X^{(\alpha)}, \mu)=\lambda_\infty (X^{(\alpha)}, \mu)$.
This together with \eqref{e:radius} yields
$ \lambda_2 (X , \mu)=\lambda_\infty (X , \mu)$ and so the conclusion of the theorem follows. \qed

\bigskip

Next, we investigate the independence of $\lambda_p(X, \mu)$
without assuming  $m(E)<\infty$.
The next result is a slight extension of \cite[Lemma 3.5]{Ta}.

\begin{lemma}\label{L:4.3}
Suppose that $\mu=\mu_1-\mu_2$, where $\mu_1$ and
$\mu_2$ are non-negative smooth measures so that
$\|G_\alpha \mu_1\|_\infty<1$ for some $\alpha >0$.
If
$$ \inf\left\{\sE (u, u) -\int_E u(x)^2 \mu (dx); \ u\in \sF
\hbox{ with } \int_E u(x)^2 m(dx)=1\right\} >0,
$$
then
$$ \inf\left\{\sE (u, u) +\int_E u(x)^2 \mu_2 (dx); \ u\in \sF
\hbox{ with } \int_E u(x)^2 \mu_1 (dx) =1 \right\} >1.
$$
\end{lemma}

\pf The proof is   the same as that for \cite[Lemma 3.5]{Ta}.
For reader's convenience, we spell out its proof here.
Let $\delta := \|G_\alpha \mu_1\|_\infty<1$.
By \eqref{eqn:en} applied to the $\alpha$-subprocess
$X^{(\alpha)}$ of $X$,  there is some $C >0$
such that
\begin{equation}\label{e:4.5}
 \int_E u(x)^2 \mu_1 (dx) \leq \delta \sE(u, u)+C  \int_E u(x)^2 m(dx)
\qquad \hbox{for } u\in \sF.
\end{equation}
Suppose that $\lambda:=\inf\left\{\sE (u, u) -\int_E u(x)^2 \mu (dx);
\ u\in \sF \hbox{ with } \int_E u(x)^2 m(dx)=1\right\} >0$; that is,
$$ \sE (u, u) -\int_E u(x)^2 \mu (dx) \geq \lambda \int_E u(x)^2 m(dx)
\qquad \hbox{for every } u\in \sF.
$$
Thus for $u\in \sF$, we have by \eqref{e:4.5} that
$$ \int_E u(x)^2 \mu_1 (dx) \leq \delta \, \sE(u, u)+\frac{C}{\lambda}
\left( \sE (u, u) -\int_E u(x)^2 (\mu_1-\mu_2) (dx)\right)
$$
and so
$$
\frac{1+(C/\lambda)}{\delta+(C/\lambda)} \int_E u(x)^2 \mu_1 (dx)
 \leq   \sE(u, u)+\frac{C/\lambda}{\delta+(C/\lambda)}
 \int_E u(x)^2  \mu_2  (dx) \leq \sE(u, u)+
 \int_E u(x)^2  \mu_2  (dx).
$$
The conclusion of the lemma now follows. \qed

\begin{remark}\label{R:4.4} \rm
 (i) Note that
$\lim_{\alpha \to \infty} \|G_\alpha \mu_1 \|_\infty=0$
for $\mu_1\in \K_\infty (X)$.
Moreover, by \cite[Proposition 2.3(ii)]{C},
\begin{equation}\label{e:4.6}
\lim_{\alpha \to \infty} \|G_\alpha \mu_1 \|_\infty<1
\qquad \hbox{for } \mu_1\in \K_1 (X).
\end{equation}

(ii) The assumption that $\|G_\alpha \mu_1\|<1$ in Lemma \ref{L:4.3}
is only used
 to deduce  inequality \eqref{e:4.5}. So the result holds under the assumption
  that   $\mu_1$ satisfies the Hardy class condition \eqref{e:4.5}
  for some $\delta <1$.
 \qed
\end{remark}

Let $\mu$ be a non-negative smooth measure and $\alpha \geq 0$.
We say $G_\alpha\mu$ is an $\alpha$-potential if $\sup_{x\in E} G_\alpha \mu (x)=0$.

\begin{lemma}\label{L:4.6}
Suppose that the process $X$ admits no killings inside, that is
$\P_x ( X_{\zeta - } \in E, \zeta <\infty)=0$ for q.e. $x\in E$.
If $\mu$ is a non-negative measure in $\K_\infty (X^{(1)})$, then
$G_\alpha \mu$ is an $\alpha$-potential for every $\alpha >0$.
\end{lemma}

\pf  Since $\mu \in \K_\infty (X^{(1)})=\K_\infty (X^{(\alpha)})$,
 for every $\eps >0$,
there is a   Borel subset $K$ and a constant $\delta >0$ so that
for every Borel subset $B\subset K$ with $\mu (B)<\eps$,
$\| G_\alpha (\1_{K^c \cup B}\mu ) \|_\infty <\eps$.
On the other hand, by \cite[Theorem 2.3.15]{CF},
there is an increasing sequence
$\{F_k, k\geq 1\}$ of closed sets so that
$\mu (\cap_{k\geq 1} F_k^c)=0$, $\mu(F_k)<\infty$,
 and $\1_{F_k}\mu$ is a measure
of finite energy for every $k\geq 1$.
Let $j\geq 1$ be large enough so that $\mu (K\setminus F_j)<\delta$.
Let $u:=G_\alpha (\1_{F_j\cap K}\mu)$. Then $u\in \sF$ and by
\cite[Corollary 3.5.3]{CF}, $\P_x ( \lim_{t\to \zeta} u(X_t)=0)=1$
for q.e. $x\in E$.  It follows that
$$
\inf_{x\in E} G_\alpha \mu (x) \leq \inf_{x\in E}
G_\alpha (\1_{ F_j \cap K}\mu ) (x) + \sup_{x\in E}
 G_\alpha (\1_{E\setminus F_j}\mu ) (x)
<\eps.
$$
Since $\eps >0$ is arbitrary, we have $\inf_{x\in E} G_\alpha \mu (x)=0$.
\qed

\begin{thm}\label{T:4.7}
 Suppose that $\mu$ is a signed smooth measure with
$\mu^+\in \K_\infty (X^{(1)})$ and $G_1\mu^-$ bounded.
\begin{description}
\item{\rm (i)} $\lambda_\infty (X, \mu)\geq \min \{\lambda_2 (X, \mu), 0\}$. Consequently, $\lambda_p (X, \mu)$ is independent of $p\in [1, \infty]$ if $\lambda_2 (X, \mu)\leq 0$.

\item{\rm (ii)} Assume in addition that  $X$ is conservative
and that either $G\mu^-$ is bounded or $G_\alpha \mu^-$ is
an $\alpha$-potential for some $\alpha>0$.
  Then $\lambda_\infty (X, \mu)=0$ if $\lambda_2 (X, \mu)>0$.
 Hence   $\lambda_p (X, \mu)$ is independent of $p\in [1, \infty]$ if and only if $\lambda_2(X, \mu)\leq 0$.
\end{description}
\end{thm}

\pf  First note that by the resolvent equation, $G_1 \mu^-$ is bounded
if and only if $G_\alpha \mu^-$ is bounded.

(i) For any $\lambda <\min\{\lambda_2 (X, \mu), 0\}$,
there is $\alpha >0$ so that $\lambda +\alpha <\min\{\lambda_2 (X, \mu), 0\}$.
Clearly,
 $$ \inf \left\{ \sE_\alpha ( u, u) -\int_E u(x)^2 ( \mu + (\lambda+\alpha)
  m) (dx); \,
 u\in \sF \hbox{ with } \int_E u(x)^2 m(dx)=1 \right\}
 = \lambda_2 (X; \mu)-\lambda  >0.
$$
Note that $(\sE_\alpha, \sF)$ is the Dirichlet form for the $\alpha$-subprocess $X^{(\alpha)}$ of $X$ and
$\|G^{(\alpha)}1\|_\infty \leq 1/\alpha$.
Since  $\mu^+ -(\mu^-+(-\lambda-\alpha) m)= \mu + (\lambda+\alpha) m$,
one has
\begin{equation}\label{e:4.8}
(\mu+(\lambda+\alpha) m)^+\leq \mu^+ \qquad \hbox{and} \qquad
(\mu+(\lambda+\alpha) m)^-\leq \mu^-+(-\lambda-\alpha) m.
\end{equation}
We thus have by Lemma \ref{L:4.3} that
$$
 \inf \left\{ \sE_\alpha ( u, u) +\int_E u(x)^2 ( (-\lambda -\alpha) m+\mu^-) (dx); \,
 u\in \sF \hbox{ with } \int_E u(x)^2 \mu^+ (dx)=1 \right\}
 >1.
$$
It follows from the elementary inequality $\frac{a}{b} \geq \frac{a+c}{b+c}$
 for $a\geq b\geq  0$ and $c>0$  that
$$
  \inf_{u\in \sF} \frac{ \sE_\alpha ( u, u) +\int_E u(x)^2
  (\mu+  (\lambda+\alpha) m)^- (dx)}{\int_E u(x)^2 (\mu+(\lambda+\alpha) m)^+ (dx)}\geq \inf_{u\in \sF} \frac{ \sE_\alpha ( u, u) +\int_E u(x)^2
(\mu^- + (-\lambda-\alpha) m) (dx)}{\int_E u(x)^2 \mu^+ (dx)}>1 .
$$
Since $\mu^+\in \K_\infty (X^{(1)})= \K_\infty (X^{(\alpha)})$
and $ G_\alpha ( (-\lambda-\alpha) m+\mu^-)$ is bounded,  it follows
from \eqref{e:4.8} and \cite[Theorem 5.4]{C}
that $(X^{(\alpha)}, \mu+(\lambda+\alpha) m)$ is
gaugeable. Let $\zeta^{(\alpha)}$ denote the lifetime of $X^{(\alpha)}$,
which can be realized as $\zeta \wedge R$ for an exponential
random variable $R$ with mean $1/\alpha$ that is independent of $X$,
and let $\E_x^{(\alpha)}$ denote the expectation under the probability law
of $X^{(\alpha)}$ starting from $x$.
Then
\begin{eqnarray*}
&&\sup_{t\geq 0} e^{\lambda  t} \|T^\mu 1\|_\infty
= \sup_{t\geq 0} \left( e^{\lambda t} \sup_{x\in E} \E_x
\left[ e^{A^\mu_t}; t<\zeta \right] \right)
= \sup_{t\geq 0} \left( e^{(\lambda +\alpha) t} \sup_{x\in E} e^{-\alpha t}\E_x \left[ e^{A^\mu_t}; t<\zeta \right] \right) \\
& = &\sup_{t\geq 0} \left( e^{(\lambda +\alpha) t} \sup_{x\in E}  \E_x \left[ e^{A^\mu_t}; t<\zeta^{(\alpha)} \right] \right)
\leq
\sup_{x\in E} \E^{(\alpha)}_x \left[ \sup_{t<\zeta} e^{(\lambda+\alpha) t +A^\mu_t} \right]<\infty,
\end{eqnarray*}
where the last inequality is due to \cite[Corollary 2.9(5)]{C}.
This implies that
$$ \lambda_\infty (X; \mu)
= -\lim_{t\to \infty} \frac1t \log \|T^\mu 1\|_\infty \geq \lambda.
$$
Since the above holds for every $\lambda <\min\{\lambda_2 (X; \mu), 0\}$,
we conclude that $\lambda_\infty (X; \mu) \geq \min\{ \lambda_2 (X; \mu), 0\}$. In particular, when $\lambda_2 (X, \mu)\leq 0$, we have
$\lambda_\infty (X, \mu)\geq \lambda_2(X, \mu)$.
This together with \eqref{e:4.2} yields that
$\lambda_p(X, \mu)$ is independent of  $p\in [1, \infty]$
when $\lambda_2(X, \mu)\leq 0$.

(ii) Since $\lambda_2 (X, \mu)>0$, we have from (i) that
$\lambda_\infty (X, \mu)\geq 0$.
Assume now that $X$ is conservative. If $G \mu^-$ is bounded,  then
$$
  \| T^\mu_t 1 \|_\infty = \sup_{x\in E} \E_x \left[ \exp ({A^\mu_t}) \right]
\geq \sup_{x\in E} \E_x \left[ \exp \left(-A^{\mu^-}_t \right) \right]
\geq \sup_{x\in E} \exp \left(-\E_x \left[ A^{\mu^-}_\infty \right]\right)
\geq  \exp \left(-\|G\mu^-\|_\infty \right) .
$$
If $G_\alpha \mu^-$ is a potential for some $\alpha >0$, then
$$
  \| T^\mu_t 1 \|_\infty
\geq \sup_{x\in E} \E_x \left[ \exp \left(-A^{\mu^-}_t \right) \right]
\geq \sup_{x\in E} \exp \left(-\E_x \left[ A^{\mu^-}_t \right]\right)
\geq  \exp \left(- \inf_{x\in E} e^{\alpha t} G_\alpha \mu^- (x) \right) =1.
$$
In either cases, we have
$$ \lambda_\infty (X, \mu)=-\lim_{t\to \infty}
 \frac1t \log \|T^\mu_t 1\|_\infty \leq 0.
$$
  Therefore  $\lambda_\infty (X, \mu)=0<\lambda_2 (X, \mu)$.
\qed

\begin{thm}\label{T:4.8}
 Suppose that $1\in \K_\infty (X^{(1)})$ and
 $\mu  \in \K_\infty (X^{(1)})$.
Then  $\lambda_p(X, \mu)$ is independent of $p\in [1, \infty]$.
\end{thm}

\pf Since $\mu \in \K_\infty (X^{(1)})$,
$\|G_1 \mu\|_\infty<\infty$.
 In view of \eqref{eqn:en}, there exists $\beta >0$ so that
 $\lambda_2 (X^{(1)}, \mu)< \beta$; or, equivalently,
$\lambda_2 (X^{(1)}, \mu +\beta m)<0$. Since $1\in \K_\infty (X^{(1)})$, we have by Theorem \ref{T:4.7} that
$\lambda_p ((X^{(1)}, \mu +\beta m)=
 \lambda_2 ((X^{(1)}, \mu +\beta m)$ for all $p\in [1, \infty]$.
 On the other hand,
 $$ \lambda_p (X^{(1)}, \mu +\beta m)=
 -\beta + \lambda_p ( X^{(1)}, \mu)
=-\beta+1 + \lambda_p (X, \mu ) \qquad \hbox{for } p\in [1, \infty],
$$
which yields the $L^p$-independence of $\lambda_p (X, \mu )$. \qed

The reason that we need to assume $1\in \K_\infty (X^{(1)})$ in Theorem
\ref{T:4.8} is because the gaugeability results for
$(X^{(1)}, \nu)$,
Theorems 2.12 and 5.2 of \cite{C},
require  $\nu^+\in \K_\infty (X^{(1)})$ and $G_1\nu^-$ bounded, and they are applied to
the measure $\nu = \mu+\beta m$. In Theorem \ref{T:4.7}, these gaugeability
results are applied to
$(X^{(\alpha)}, \nu)$ for measure
 $\nu = \mu +(\lambda + \alpha) m$ with
$\lambda +\alpha <0$ so we do not need to assume $1\in \K_\infty (X)$.

\medskip

The following sufficient condition for $1\in \K_\infty (X^{(1)})$ is established in \cite[Theorem 4.2]{C} (together with its proof).

\begin{lemma}\label{L:4.9}
Suppose that $\alpha>0$ and that $G_{\alpha}$ maps bounded functions into continuous functions. If for every $\eps >0$, there is a compact set $K  \subset E$
such that $\sup_{x\in E} G_{\alpha} \1_{K^c}(x) <\eps$, then
$1\in \K_\infty (X^{(\alpha)})=\K_\infty (X^{(1)})$.
\end{lemma}

\begin{remark}\label{R:4.10} \rm
(i) Theorem \ref{T:4.7}(ii) extends the main result (Theorem 3.1) of \cite{T},
where it is shown by a large deviation argument that,
for any $m$-symmetric irreducible conservative
Feller process $X$ that has jointly continuous transition density function
and signed measure $\mu \in \K_\infty (X)$, $\lambda_p(X, \mu)$ is
independent of $p\in [1, \infty]$ if and only if $\lambda_2 (X, p)\leq 0$.
Theorem \ref{T:4.8} extends a corresponding result in \cite{T2}.
Our approach also reveals where the role of conservativeness of $Y$
is played in part (ii) of Theorem \ref{T:4.7} in
connection with  part (i). See \cite{C, Ta} for related results
on the $L^p$-independence of the spectral radius of the transition
semigroup of $X$ (that is, for $\lambda_p(X, 0)$
corresponding to $\mu=0$).

(ii) Assume that $X$ is an $m$-symmetric irreducible process $X$ satisfying strong Feller property (that is, its transition semigroup maps bounded Borel
measurable functions into bounded continuous functions) and
the following tightness assumption:
for every $\eps>0$, there is a compact subset $K$ so that
$\sup_{x\in E} G_1 \1_K(x)\leq \eps$.
For such a process, as an application of
a large deviation result established in
\cite[Theorem 1.1]{T2}, it is shown in \cite{T2} that
$\lambda_p(X, \mu)$ is independent of $p\in [1, \infty]$ for every
$\mu \in \K_\infty (X)$. This result is a special case of our Theorem
\ref{T:4.8} in view of Lemma \ref{L:4.9}.

(iii) By the same argument as that for \cite[Proposition 4.1]{T2}
that for smooth measure $\mu$ with $\mu^+\in \K_1(X)$ and
$G \mu^+$ bounded,
$$ \liminf_{t\to \infty} \frac1t \log \E_x
\left[ e^{A^\mu_t}; t<\zeta \right]
\geq -\lambda_2 (X, \mu), \qquad x\in E.
$$
So whenever $\lambda_2 (X, \mu)=\lambda_\infty (X, \mu)$, one has
for every $x\in E$,
$$
 \lim _{t\to \infty} \frac1t \log \E_x
 \left[ e^{A^\mu_t}; t<\zeta \right]
=\lim _{t\to \infty} \frac1t \log \sup_{x\in E}
\E_x \left[ e^{A^\mu_t}; t<\zeta \right] = -\lambda_2 (X, \mu).
$$
\qed
\end{remark}

\section{Spectral bounds for non-local Feynman-Kac semigroups}\label{S:5}

Throughout this section,
 $F$ is a bounded symmetric function  in the Kato class ${\bf J} (X)$.
 Let $M$ be the purely discontinuous square-integrable martingale additive functional
of $X$ with
$$
 \Delta M_t=e^{F(X_{t-}, X_t)}-1,
$$
 which has the expression \eqref{e:3.2}.
Let $\Exp (M)$ be the Dol\'eans-Dade exponential martingale of $M$.
It defines a family of probability measures $\{\bQ_x, x\in E\}$ by
$d\bQ_x/ d\bP_x = \Exp (M)_t $ on $\sF_t$.
For emphasis,
the Girsanov transformed process $\{X_t, \bQ_x\}$ is denoted by $Y$.
The process $Y$ is still $m$-symmetric and  its associated
Dirichlet form on $L^2(E; m)$ is $(\sE^Y, \sF)$, where
\begin{equation}\label{e:5.1}
 \sE^Y (u, u) = \sE (u, u)+\frac12 \int_{E\times E} (u(x)-u(y))^2
 \left(e^{F(x, y)}-1 \right) N(x, dy)\mu_H (dx)
\end{equation}
(see \cite{CS}). So the symmetric process $Y$ has L\'evy system $\left( e^{F(x,
y)}N(x, dy), H \right)$.
If $F\in \J_\infty (X)$, then by Remark \ref{R:3.4}(ii), $\Exp (M)$ is a uniformly integrable martingale
under $\P_x$ for every $x\in E$.
\medskip

\begin{lemma}\label{L:5.1}
 If $A$ is a PCAF of $X$ with Revuz measure $\nu$, then
the Revuz measure of $A$ as a PCAF of $Y$ is still $\nu$.
\end{lemma}

\pf The proof is exactly the same as that for Lemma 4.4 in \cite{Fi}
so it is omitted here. \qed

The next result says that
the corresponding Kato classes become larger after Girsanov transform.

\begin{thm}\label{T:5.2} Let
$F\in {\bf J}_\infty (X)$ be symmetric
and $Y$ be the above Girsanov transformed process
in terms of $F$. Then
$${\bf K}_\infty (X)\subset {\bf K}_\infty (Y) \qquad
\hbox{and} \qquad  {\bf J}_\infty (X)\subset {\bf J}_\infty (Y).
$$
Moreover, if $\nu$ is a non-negative measure so that $G\nu$ is bounded,
then so is $G^Y\nu$.
\end{thm}

\pf For notational convenience, let $Z_t=\Exp (M)_t$.
Since $\bE_x\left[ Z_\zeta \right] =\bE_x [ Z_0]= 1$ for every
$x\in E$, $(X, F- A^F)$ is gaugeable. By the Super Gauge theorem
(Theorem \ref{T:2.3}(ii)), there is an $\eps>0$ so that
 \bee\label{e:superg}
 c_0:=\sup_{x\in E} \left( \E_x \left[ Z_\zeta^{1+\eps}\right]\right)^{1/(1+\eps)}
  <\infty.
 \eee
 Clearly, $Y$ has a Green function $G^Y(x, y)$ defined by
 $$ \int_E G^Y(x, y) f(y)=\E_x \left[ \int_0^\infty f(Y_s) ds\right]
 =\E_x \left[ Z_\zeta \int_0^\zeta f(X_s) ds\right].
 $$
 In fact, $G^Y(x, y)=\E_x^{y} [ Z_\zeta] G(x, y)$, where
 $\E_x^y$ is the expectation under the law of $\P_x^y$ which
 is obtained from $\P_x$ through Doob's $h$-transform with
 $h(z)=G(z, y)$; see \cite{C}. Hence for each fixed $y$,
 $x\mapsto G^Y(x, y)$ is an excessive function of $Y$.

Let $k\geq 2$ be an integer so that $p:=k/(k-1)<1+\eps$. Then for
any positive smooth measure $\nu$ with $\|G\nu \|_\infty = \|
\E_{\cdot} A^\nu_{\zeta}\|_\infty <\infty$, by H\"older's inequality
and \eqref{e:superg},
\begin{equation}\label{e:kato2}
G^Y \nu (x) = \E_x \left[ Z_\zeta A^\nu_{\zeta}\right] \leq \left(
\E_x \left[ Z^p_\zeta\right] \right)^{1/p} \left( \E_x \left[
(A^\nu_{\zeta} )^k \right] \right)^{1/k} \leq c_0  \left( k!
\right)^{1/k} \, \| G\nu\|_\infty.
\end{equation}
This implies  that ${\bf K}_\infty (X)\subset {\bf K}_\infty (Y)$
and so  ${\bf J}_\infty (X)\subset {\bf J}_\infty (Y)$.
In particular, \eqref{e:kato2} implies that $G^Y\nu$ is bounded if
$G\nu$ is bounded.
\qed

\bigskip

Clearly, the 1-subprocess $Y^{(1)}$ of $Y$ can be obtained
from $X^{(1)}$ through the Girsanov transform $\Exp (M)$.
When $F\in \J_\infty (X^{(1)})$, we have by applying Theorem \ref{T:5.2}
to  $X^{(1)}$   that
\begin{equation}\label{e:5.4a}
{\bf K}_\infty (X^{(1)})\subset {\bf K}_\infty (Y^{(1)}) \qquad
\hbox{and} \qquad  {\bf J}_\infty (X^{(1)})\subset {\bf J}_\infty (Y^{(1)}).
\end{equation}

Assume that $\mu$ is a signed smooth measure with $\mu^+\in
\K  (X)$ and $G\mu^-$ bounded, and $F\in \J  (X)$ symmetric.
Define the non-local Feynman-Kac semigroup
$$T^{\mu, F}_tf (x):=\bE_x \Big[ \exp\Big( A^\mu_t+\sum_{0<s\leq t} F(X_{s-}, X_s) \Big) f(X_t)\Big], \qquad t\geq 0.
$$
It follows from the proof of \cite[Proposition 2.3]{CS1}
and H\"older  inequality that
$\{T^{\mu, F}_t; t\geq 0\}$ is a strongly continuous semigroup in
  $L^p(E; m) $ for every $1\leq p\leq \infty$. Moveover, it is easy to verify that
$T^{\mu, F}_t$ is a symmetric operator in $L^2(E; m)$.
The $L^p$-spectral bound of $\{ T^{\mu, F}_t; t\geq 0\}$ is defined to
be
$$ \lambda_p(X,  \mu+F):= -\lim_{t\to \infty}
    \frac1t \log \| T^{\mu, F}_t\|_{p, p}.
$$

The purpose of this section is to give  necessary and sufficient conditions
for $\lambda_p(X,  \mu+F)$ to be independent of $1\leq p\leq \infty$.

Assume from now that $F$ is a bounded symmetric function in
 $\J_\infty (X^{(1)})$.
Let $A^F$ be the continuous additive functional defined by \eqref{e:3.3},
whose Revuz measure is
$$
\nu_F (dx):= \left(\int_E (e^{F(x, y)}-1) N(x, fy) \right) \mu_H (dx).
$$
By \eqref{e:3.4}, we have
 \begin{equation}\label{e:5.4}
 T^{\mu, F}_t f(x)= \E_x \left[ \Exp (M)_t \, e^{A^\mu_t+A^F_t} \, f(X_t)\right] = \E_x^\bQ
 \left[ e^{A^\mu_t+A^F_t} \, f(Y_t) \right] =: Q_t f(x).
\end{equation}
Thus for any $1\leq p\leq \infty$ and $t\geq 0$,
 \bee\label{e:5.5}
 \| P_t^{\mu, F}\|_{p,p} = \|Q_t\|_{p, p} \qquad \hbox{ and so } \qquad
\lambda_p(X;  \mu+F)=\lambda_p(Y;  \mu+\nu_{F}).
 \eee
  Since $F$ is bounded, $c_1 |F| \leq |e^F-1|\leq c_2 |F|$.
Thus in view of Lemma \ref{L:5.1} and \eqref{e:5.4a},
 \bee\label{e:5.6}
 \hbox{the signed Revuz measure $\nu_F$ of } A^F \hbox{ belongs to  }
  {\bf K}_\infty (X^{(1)})\subset {\bf K}_\infty (Y^{(1)}).
 \eee
In particular, it follows from \eqref{e:4.1}, \eqref{e:5.1}
and \eqref{e:5.5}-\eqref{e:5.6} that
\begin{eqnarray}
 && \lambda_2 (X;  \mu+F)=\lambda_2(Y;  \mu+\nu_F ) \nonumber \\
&=& \inf \left\{ \sE^Y(u, u)-\int_E u(x)^2 \left(\int_E
\left( e^{F(x, y)}-1 \right)N(x, dy) \right) \mu_H(dx)
 -\int_E u(x)^2 \mu (dx); \right. \nonumber \\
&& \hskip 0.5truein \left.  u\in \sF \hbox{ with } \int_E u(x)^2 m(dx)=1 \right\} \nonumber \\
&=& \inf \left\{ \sE (u, u) -
 \int_{E\times E} u(x)u(y) \left( e^{F(x, y)}-1\right) N(x, dy) \mu_H(dx)-\int_E u(x)^2 \mu (dx); \right. \nonumber \\
&& \hskip 0.5truein \left.  u\in \sF \hbox{ with } \int_E u(x)^2 m(dx)=1 \right\}.  \label{e:5.7}
\end{eqnarray}

\bigskip

We start with an analogy of Corollary \ref{C:radius}.

\begin{thm}\label{T:5.3}
Assume that \eqref{e:lp3} holds and  $m(E)<\infty$.
Let $\mu$ be  a signed smooth measure
with $\mu^+\in  \K_\infty (X^{(\alpha)})$
and $G_{\alpha } \mu^-$ bounded for some $\alpha \geq 0$,
and $F\in \J_\infty (X^{(\alpha)})$ symmetric.
Then
$\lambda_p(X, \mu+F)$ is independent of $p\in [1, \infty]$.
\end{thm}

\pf  By the same reasoning as that in Remark \ref{R:4.2},
 we have
 \begin{equation}\label{e:5.9}
 \lambda_p (X^{(\alpha)}, \mu+F) = \lambda_p(X, \mu+F)+ \alpha
  \qquad \hbox{for every }p\in [1, \infty].
 \end{equation}
So  it suffices to show that $\lambda_p (X^{(\alpha)}, \mu+F)$
is independent of $p\in [1, \infty]$. By regarding $X^{(\alpha)}$
as $X$, we may assume, without loss of generality, that
the condition of the theorem holds with $\alpha =0$.

Since $P_{t_0}$ is a bounded linear operator from $L^2(E; m)$
  to $L^\infty (E; m)$, by duality, $P_{t_0}$ is a bounded linear operator
  from $L^1(E; m)$ to $L^2(E; m)$. Hence $P_{2t_0}: L^1(E; m)\to L^\infty (E; m)$ is bounded. Let
  $$
  Z_t=\Exp (M)_t=\exp \left( \sum_{0<s\leq t} F(X_{s-}, X_s)-A^F_t\right),
  \qquad t\geq 0 .
  $$
Since $F\in \J_\infty (X)$, it follows from  Khasminskii's inequality
and the Markov property that (cf. \cite[(3.11)]{C2}) that there are constants
$c_1, c_2>0$ so that
$\sup_{x\in E} \E \left[ Z_t^2\right] \leq c_1e^{c_2t}$ for every $t>0$.
Denote by $Y$ the Girsanov transformed process of $X$ via $Z$.
 Then for every $f\in L^2(E; m)$,
$$ |P^Y_{2t_0} f(x)|:= |\E_x [ f(Y_{2t_0})|=
 \E_x \left[ M_{2t_0} f(X_t)\right]\leq \left( \E_x [M^2_{2t_0}]
 \, \E_x \left[ f(X_{2t_0})^2\right]\right)^{1/2}
 \leq c\, \|f\|_{L^2(E; m)}.
 $$
 This proves that condition \eqref{e:lp3} holds for $Y$
  with $2t_0$ in place of $t_0$.
 Since
 $$(\mu+\nu_F)^+\leq \mu^+ + (\nu_F)^+ \quad \hbox{ and } \quad
 (\mu+\nu_F)^- \leq \nu^- + (\nu_F)^-,
 $$
 we deduce by  Theorem \ref{T:5.2} and \eqref{e:kato2} that
 $(\mu+\nu_F)^+\in \K_\infty (Y)$ and $G^Y (\mu+\nu_F)^-\, $ is bounded.
 Hence by \eqref{e:5.5} and Corollary \ref{C:radius},
  $\lambda_p(X, \mu+F)=\lambda_p (Y, \mu+\nu_F)$ is independent
  in $p\in [1, \infty]$. \qed

The next result is a non-local Feynman-Kac semigroup counterpart
of Theorem \ref{T:4.7}.

\begin{thm}\label{T:5.4}
 Suppose that $\mu$ is a signed smooth measure with
$\mu^+\in \K_\infty (X^{(1)})$ and $G_1\mu^-$ bounded,
and $F\in \J_\infty (X^{(1)})$ symmetric.
\begin{description}
\item{\rm (i)} $\lambda_\infty (X,  \mu+F)\geq \min \{\lambda_2 (X, \mu+F), 0\}$. Consequently, $\lambda_p (X,  \mu+F)$ is independent of $p\in [1, \infty]$ if $\lambda_2 (X, \mu+F)\leq 0$.
\item{\rm (iii)} Assume in addition that  $X$ is conservative
and that   $\mu\in \K_\infty (X^{(1)})$.
 Then $\lambda_\infty (X, \mu +F )=0$
 if   $\lambda_2 (X, \mu+F)>0$.
 Hence   $\lambda_p (X, \mu +F)$ is independent of $p\in [1, \infty]$ if and only if $\lambda_2(X, \mu +F )\leq 0$.
\end{description}
\end{thm}

\bigskip

\pf   For notational convenience, let $Z_t:=\Exp (M)_t$.
By Remark \ref{R:3.4}(ii),  $\{Z_t, t\geq 0\}$ is a uniformly integrable martingale under each $\bP_x$. It follows that the Girsanov transformed process $Y$ is transient and has a Green function
$G^Y$. Furthermore, $Y$ is conservative if so is $X$.
  It is clear that $Y$ is $m$-irreducible since $\Exp(M)_t>0$ a.s..
 Note that it follows from Theorem \ref{T:5.2} applied to
 the subprocess $X^{(1)}$ that if $G_1\mu^-$ is bounded,
 then so is $G^Y_1 \mu^-$. Thus in view of
  \eqref{e:5.4a} and Lemma \ref{L:4.6},
  $\mu+\nu_F$ satisfies the condition
  of Theorem \ref{T:4.7} for the symmetric process $Y$.
   The conclusion of the theorem now follows from
Theorem \ref{T:4.7} applied to $(Y, \mu+\nu_F)$
and relation \eqref{e:5.5}. \qed

The following theorem extends Theorem \ref{T:4.8} to non-local
Feynman-Kac semigroups.

\begin{thm}\label{T:5.5}
 Suppose that $1\in \K_\infty (X^{(1)})$,
 $\mu  \in \K_\infty (X^{(1)})$ and $F\in \J_\infty (X^{(1)})$
 symmetric.
Then  $\lambda_p(X, \mu+F)$ is independent of $p\in [1, \infty]$.
\end{thm}

\pf Let $Y^{(1)}$ be the Girsanov transformed process from $X^{(1)}$
via the function $F$. By the same reason as that for Theorem \ref{T:5.4},
we can apply Theorem \ref{T:4.8} to $(Y^{(1)}, \mu+\nu_F)$ to conclude
that $\lambda_p (Y^{(1)}, \mu+\nu_F)$ is independent of $p\in [1, \infty]$. Consequently, in view of \eqref{e:5.5} applied to $X^{(\alpha)}$,
$\lambda_p (X^{(1)}, \mu+ F)$ is independent of $p\in [1, \infty]$.
The conclusion of the theorem follows once one notices that
$\lambda_p (X^{(1)}, \mu+ F)= 1 +  \lambda_p (X , \mu+ F)$
for every $p\in [1, \infty]$.
\qed

\begin{remark}\label{R:5.6}
 \rm (i) When $\mu =0$, the conclusion of Theorem \ref{T:5.4}(i)
 recovers and extends the main result of \cite{Taw}, which was established by
   using a large deviation approach. The latter extends
   an earlier result of \cite{TT} where $X$ is a rotationally
   symmetric $\alpha$-stable process.

(ii) It follows from \eqref{e:5.4} and Remark \ref{R:4.10}(iii)
that for any  $\mu \in \K_\infty (X)$ and  symmetric
$F\in \J_\infty (X)$,
$$ \liminf_{t\to \infty} \frac1t \log \E_x \left[ \exp
\left( A^\mu_t+\sum_{0<s\leq t} F(X_{s-}, X_s) \right);
t<\zeta \right] \geq -\lambda_2 (X, \mu+F), \qquad x\in E.
$$
Thus whenever $\lambda_2 (X, \mu +F)=\lambda_\infty (X, \mu+F)$,
we have for every $x\in E$,
$$ \lim _{t\to \infty} \frac1t \log \E_x \left[ \exp \left(
A^\mu_t+\sum_{0<s\leq t} F(X_{s-}, X_s) \right); t<\zeta \right] = -\lambda_2 (X, \mu+F).
$$

(iii) The idea of using pure Girsanov transform to reduce a
non-local Feynman-Kac
transform to a continuous (local) Feynman-Kac transform of a new process
is applicable in many other situations. For example, using this idea,
one can easily deduce a large deviation result for general
non-local Feynman-Kac functionals \cite[Theorem 2.1]{TT3}
from the corresponding result of  local Feynman-Kac transforms
\cite[Theorem 1.1]{T2}.
In fact, such an approach shows that the large deviation result
in \cite[Theorem 2.1]{TT3} holds in fact for $\mu\in \K_\infty (X)$
and symmetric $F\in \J_\infty (X)$,
while \cite[Theorem 2.1]{TT3} requires
$\mu\in \K_\infty(X)$ and symmetric $F\in {\bf A}_2(X)$,
a subclass of $ \J_\infty (X)$
introduced in \cite{C}.
It should be mentioned that \cite[Theorem 2.1]{TT3}
implies  the independence of $\lambda_p(X, \mu+F)$ in $p\in [1, \infty]$
under the same condition on $X$ as in \cite{T2}
(see Remark \ref{R:4.10}(ii) above) with
$\mu \in \K_\infty (X)$ and symmetric $F\in {\bf A}_2(X)$.
\qed
\end{remark}

We refer the reader to \cite[Section 5]{C3} for  concrete examples
of functions in Kato classes
$\K_\infty (X^{(1)})$ and $\J_\infty (X^{(1)})$.
For example, it is shown in \cite{C3} that when $X$ is a
symmetric $\alpha$-stable-like process on a global $d$-set $E$
with $d$-measure $m$, then
$L^p(E; m)\subset  \K_\infty (X^{(1)})$ for every $p>d/\alpha$
when $\alpha \leq d$, and
$L^p(E; m)\subset \K_\infty (X^{(1)})$ for every $p\geq 1$
when $0<d<\alpha$. It is further shown there that when
$X$ is a symmetric diffusion on $\bR^d$ associated with a
uniformly elliptic and bounded divergence form operator,
then  $ L^p(\bR^n; dx)\subset \K_\infty (X^{(1)})$
for every $p>n/2$ when $n\geq 3$, and
$L^p(\bR^n; dx)\subset \K_\infty (X^{(1)})$
for every $p\geq 1$ when $d=1$ or $2$.

\bigskip

\nin{\bf Acknowledgement.} The author thanks M. Takeda for helpful comments
on an earlier version of this paper.

\vskip 0.3truein

\begin{singlespace}

\end{singlespace}

\vskip 0.3truein

{\bf Zhen-Qing Chen}

Department of Mathematics, University of Washington, Seattle,
WA 98195, USA

E-mail: \texttt{zqchen@uw.edu}

\end{document}